\documentclass[11pt,twoside,draft]{article}
\usepackage{amsmath}
\usepackage{amssymb}
\usepackage{graphicx}
\usepackage{amsfonts}
\usepackage{setspace}
\usepackage[figuresright]{rotating}
\usepackage{lscape}
\usepackage{longtable}
\usepackage{slashbox}

\newtheorem{thm}{Theorem}

\newtheorem{lem}[thm]{Lemma}

\newcommand{\ncm}{\newcommand}
\ncm{\bc}{\begin{center}} \ncm{\ec}{\end{center}}
\ncm{\be}{\begin{equation}} \ncm{\ee}{\end{equation}}
\ncm{\ba}{\begin{array}} \ncm{\ea}{\end{array}}
\ncm{\beq}{\begin{eqnarray}} \ncm{\eeq}{\end{eqnarray}}
\ncm{\lab}[1]{\label{#1}} \ncm{\re}[1]{(\ref{#1})}
\ncm{\nn}{\nonumber}

\def\ds{\displaystyle}

\def\a{\alpha}
\def\b{\beta}
\def\ga{\gamma}
\def\de{\delta}
\def\s{\sigma}
\def\vp{\varphi}
\def\ve{\varepsilon}

\begin{document}
\large

\begin{center}
{ \Large \bf Estimation of errors of quadrature formula for singular integrals of Cauchy type with special forms} \\
\end{center}

\begin{center}

\begin{tabular}{|c|}
  \hline
  Israilov Maruf Israilovich \\
  \hline
\end{tabular}

\begin{center}
{\small Institute of Mathematics, Academy of Science of Uzbekistan.}
\end{center}

{\small \textsl{Translator: Associate Prof. Dr. Zainidin K. Eshkuvatov, \\ "Department of Mathematics, Faculty of Science, Universiti Putra Malaysia"}
\footnote[1]{This paper is dedicated to the memory of Professor Israilov M.I. (1934-2010). It was published in collection of papers entitled "Differential Equations and Inverse Problems, Press FAN, Tashkent, 1986, pp. 236-258 (in Russian)}}.

\end{center}



\begin{abstract} In this work, we consider the singular integrals of Cauchy type of the forms
\begin{equation} \label{eq1}
 \ds J(f,x)= \frac{\sqrt{1-x^2}}{\pi}\int_{-1}^1\frac{f(t)}{\sqrt{1-t^2}(t-x)}\,dt, \ \ \qquad-1<x<1.
\end{equation}
and
\begin{equation} \label{eq2}
\ds \Phi(f,z)=
-\frac{\sqrt{z^2-1}}{\pi}\int_{-1}^1\frac{f(t)}{\sqrt{1-t^2}(t-z)}\,dt,
\ \ \qquad z\notin [-1,1].
\end{equation}
which are understood as Cauchy principal value integrals. Quadrature formulas (QFs) for singular integrals (SIs) \re{eq1}
and \re{eq2} are of the forms
\begin{equation} \label{eq3}
 \ds J(f,x)= \sum_{k=0}^{N}A_k(x)f(t_k)+ R_N(f,x), \ \ \qquad-1<x<1.
\end{equation}
and
\begin{equation} \label{eq4}
\ds \Phi(f,z)= \sum_{k=0}^{N}B_k(z)f(t_k)+ R_N^*(f,z), \ \qquad
z\notin [-1,1].
\end{equation}
where $z$ is complex variable with $|Re(z)|>1$. With the help of linear spline interpolation, we have proved the rate of convergence
of the errors of QFs \re{eq3} and \re{eq4} for different classes (i.e. $H^\a([-1,1],K), C^{m,\a}[-1,1], W^r[-1,1]$) of density
function $f(t)$. It is shown that approximation by spline possesses more advantages than other kinds of approximation: it requires the minimum smoothness of density function $f(x)$ to get good order of decreasing errors.

\end{abstract}

\section{Introduction}
The importance of singular integrals (SIs) of the form \re{eq1} and \re{eq2} and their numerical solution  are given in many researchers work (\cite{BL}-\cite{Ga2}) and literatures cited therein. Many of them are based on the approximation of density function $f(t)$ with Chebyshev polynomials.

Note that in \re{eq2}, the function $\sqrt{z^2-1}$  is understood as a single-valued branch in the plane of complex variable with cut along the interval [-1,1] such that $\sqrt{z^2-1}= z+O(z^{-1})$ for the large $z>0$. In the future,  under $W=arcsinz$ only a branch of the function for which $|Re(W)|< \frac{\pi}{2}$ will be understood.

In this paper, we construct efficient quadrature formulas (QFs) for SIs \re{eq1} and \re{eq2} using linear spline interpolation. Obtained QFs provide uniform convergence for any singular point $x\in(-1,1)$ and any $z \notin [-1,1]$.

\section{Construction of the quadrature formula}
In order to write the exact form of coefficients of the quadrature formula \re{eq3} and \re{eq4}, we introduce the following notations
$$
G_k(x)=\left|\frac{t_k\sqrt{1-x^2}-x\sqrt{1-t_k^2}}{\sqrt{1-x^2}+\sqrt{1-t_k^2}}\right|,
$$
$$
g_k=\frac{1}{\pi h_k}(arcsin\,t_{k+1}-arcsin\,t_{k}),
$$
$$
F_k(z)=\frac{1}{\pi h_k}
\left(arcsin\frac{zt_{k+1}-1}{z-t_{k+1}}-arcsin\frac{zt_{k}-1}{z-t_{k}}\right).
$$
If $x\not=t_k$, then the coefficients of QFs \re{eq3} are computed by the formulas
\begin{eqnarray}
\left.
\begin{array}{lll}
&& \ds A_0(x)= \frac{t_1-x}{\pi h_0}lnG_1(x)-\sqrt{1-x^2}g_0, \\[5mm]
&& \ds A_k(x) = \ds \frac{t_{k+1}-x}{\pi h_k}
ln\frac{G_{k+1}(x)}{G_{k}(x)} + \frac{x-t_{k-1}}{\pi h_{k-1}}
ln\frac{G_{k}(x)}{G_{k-1}(x)} \\[5mm]
&& \ds \qquad \qquad - \sqrt{1-x^2}(g_k-g_{k-1}), \ \ k=1,...,N-1 \\[5mm]
&& \ds A_{N}(x) = \frac{t_{N-1}-x}{\pi
h_{N-1}}ln\,G_{N-1}(x)+\sqrt{1-x^2}g_{N-1}. \label{eq5}
\end{array}
\right\}.
\end{eqnarray}
As $G_k(\pm 1)=1$ for all $k$, then $A_k(\pm1)=0 $ for $k=0,...,N$. These correspond with the fact in \cite{Pi} that
$J(f,x)|_{x=\pm 1}=0$ is independent from the value of $f(\pm 1)$ .

If $x$ coincides with the nodes $t_k, (k=1,...,N-1)$, then the coefficients $A_j(t_k), j\not= k-1,k,k+1$ are computed by
\re{eq5}. If $k$ in \re{eq5} is replaced by $k-1$ and $k+1$ and $x=t_k$ is put, then coefficients $A_{k-1}(t_k)$ and $A_{k+1}(t_k)$
are again computed respectively by \re{eq5} and for $A_k(t_k)$ we have
\be \label{eq6}
A_k(t_k)=\frac{1}{\pi}ln\frac{G_{k+1}(t_k)}{G_{k-1}(t_k)} -
\sqrt{1-t_k^2}(g_k-g_{k-1}), \ \ k=1,...,N-1.
\ee
Coefficients of the QFs \re{eq4} have the form
\begin{eqnarray}
\left.
\begin{array}{lll}
&& \ds B_0(z)= (z-t_1)F_0(z)+\sqrt{z^2-1}g_0, \\[5mm]
&& \ds B_k(z) = \ds (z-t_{k+1})F_k(z) -(z-t_{k-1})F_{k-1}(z)
\\[5mm]
&& \ds \qquad \qquad +  \sqrt{z^2-1}(g_k-g_{k-1}), \ \ k=1,...,N-1 \\[5mm]
&& \ds B_{N}(z) = (t_{N-1}-z)F_{N-1}(z) - \sqrt{z^2-1}g_{N-1}.
\label{eq7}
\end{array}
\right\}.
\end{eqnarray}
Let us derive coefficients of QFs which are given by \re{eq5} and \re{eq7}. As we know the linear spline $S_N(t)$ interpolating the given function $f$ on the grid $\Delta: \ \ -1=t_0< t_1< ... < t_{N-1}< t_N=1$ for $t \in [t_j,t_{j+1}]$ has the form
\begin{eqnarray}
S_N(t)=\frac{1}{h_k}\Bigl[(t_{k+1}-t)f(t_k)+(t-t_k)f(t_{k+1})\Bigr] \label{eq8}
\end{eqnarray}
Replacing $f(t)$ in \re{eq1} with $S_N(t)$ we have
\begin{eqnarray}
&& \ds J(S_N,x)=\frac{\sqrt{1-x^2}}{\pi}\sum\limits_{k=0}^{N-1}\frac{1}{h_k}\int\limits_{t_k}^{t_{k+1}}
\frac{\Bigl[(t_{k+1}-t)f(t_k)+(t-t_k)f(t_{k+1})\Bigr]}{\sqrt{1-t^2}(t-x)}dt. \nn \\
&& \qquad \quad \ds =\frac{\sqrt{1-x^2}}{\pi}
\left[\frac{1}{h_0}\int\limits_{-1}^{t_{1}}
\frac{(t_{1}-t)dt}{\sqrt{1-t^2}(t-x)}f(t_0) + \frac{1}{h_{N-1}}\int\limits_{t_{N-1}}^{1}
\frac{(t-t_{N-1})dt}{\sqrt{1-t^2}(t-x)}f(t_N) \right. \nn \\
&& \qquad \quad \ds + \sum\limits_{k=1}^{N-2}\frac{1}{h_k}\left(\int\limits_{t_k}^{t_{k+1}}
\frac{(t-t_k)dt}{\sqrt{1-t^2}(t-x)}+\int\limits_{t_{k+1}}^{t_{k+2}}
\frac{(t_{k+2}-t)dt}{\sqrt{1-t^2}(t-x)} \right)f(t_{k+1}) . \label{eq9}
\end{eqnarray}
Introducing notations
$$
J(k,x)=\int\limits_{t_k}^{t_{k+1}}\frac{dt}{\sqrt{1-t^2}(t-x)}, \ \
J_1(k,x)=\int\limits_{t_k}^{t_{k+1}}\frac{tdt}{\sqrt{1-t^2}(t-x)},
$$
and using easy checking formulas
$$
\int\frac{xdt}{\sqrt{1-t^2}(\sqrt{1-x^2}+\sqrt{1-t^2})}=
ln\frac{1+\sqrt{1-x^2}\sqrt{1-t^2}+xt}{\sqrt{1-x^2}+\sqrt{1-t^2}} +C,
$$
$$
\int\frac{tdt}{\sqrt{1-t^2}(\sqrt{1-x^2}+\sqrt{1-t^2})}=
-ln(\sqrt{1-x^2}+\sqrt{1-t^2})+C,
$$
obviously we have
$$
J(k,x)=\frac{1}{1-x^2}ln\left|\frac{(t-x)(1+\sqrt{1-x^2}\sqrt{1-t^2}+xt)}
{(\sqrt{1-x^2}+\sqrt{1-t^2})^2}\right|_{t=t_k}^{t=t_{k+1}},
$$
writing
$$
(t-x)(1+\sqrt{1-x^2}\sqrt{1-t^2}+xt)=(\sqrt{1-x^2}+\sqrt{1-t^2})(t\sqrt{1-x^2}-x\sqrt{1-t^2}),
$$
we obtain
\begin{eqnarray}
J(k,x)=\frac{1}{1-x^2}ln\left|\frac{t\sqrt{1-x^2}-x\sqrt{1-t^2}}
{(\sqrt{1-x^2}+\sqrt{1-t^2})^2}\right|_{t=t_k}^{t=t_{k+1}}=
\frac{1}{1-x^2}ln\frac{G_{k+1}}{G_k(x)}. \label{eq10}
\end{eqnarray}
Next, having the following relation
\begin{eqnarray}
J_1(k,x)=xJ(k,x)+\pi h_kg_k. \label{eq11}
\end{eqnarray}
Eq. \re{eq9} can be rewritten as
\begin{eqnarray}
&& \ds J(S_N,x)=\frac{\sqrt{1-x^2}}{\pi}\left\{\frac{1}{h_0}\Bigl[t_1J(0,x)-J_1(0,x)\Bigr]f(-1) \right. \nn\\
&& \ds \qquad \qquad +\frac{1}{h_{N-1}}\Bigl[J_1(N-1,x)-t_{k-1}J_1(N-1,x)\Bigr]f(1)\nn \\
&& \ds \qquad \qquad +\sum\limits_{k=1}^{N-1}\left[\frac{1}{h_k}\Bigl(t_{k+1}J(k,x)-J_1(k-1,x)\Bigr) \right.\nn \\
&& \ds \qquad \qquad + \left. \frac{1}{h_{k-1}}\Bigl(J_1(k-1,x)-t_{k-1}J(k-1,x)\Bigr)\right]f(t_k) \label{eq12}
\end{eqnarray}
Substituting \re{eq10} and \re{eq11} into \re{eq12} and simplifying the expressions, we arrive at \re{eq5} for finding the coefficients $A_k(x)$ 
of the QFs \re{eq3}. Furthermore, we can derive the Eq. \re{eq6} from \re{eq12} and \re{eq10}-\re{eq11} as follows
\begin{eqnarray}
&& \ds A_k(t_k)=\frac{\sqrt{1-x^2}}{\pi}\left\{\frac{1}{h_k}\Bigl(t_{k+1}J(k,t_{k+1})-J_1(k,t_k)\Bigr) \right.\nn \\
&& \qquad \qquad \left. +\frac{1}{h_{k-1}}\Bigl(J_1(k-1,t_{k})-t_{k-1}J(k-1,t_k)\Bigr) \right\} \nn \\
&& \qquad \quad  = \frac{\sqrt{1-x^2}}{\pi}\left[J(k,t_{k})-J(k-1,t_k)-\pi (g_k-g_{k-1})\right] \nn \\
&& \qquad \quad  = \frac{1}{\pi}ln\frac{G_{k+1}(t_k)}{G_{k-1}(t_k)} -
\sqrt{1-t_k^2}(g_k-g_{k-1}), \ \ k=1,...,N-1. \label{eq13}
\end{eqnarray}
In order to derive the coefficients $B_k$ of the QFs \re{eq4} we use the combination of the integrals
\begin{eqnarray}
J^*(k,z)=\int\limits_{t_k}^{t_{k+1}}\frac{dt}{\sqrt{1-t^2}(t-z)}, \ \
J_1^*(k,z)=\int\limits_{t_k}^{t_{k+1}}\frac{t \,dt}{\sqrt{1-t^2}(t-z)}. \label{eq14}
\end{eqnarray}
First $J^*(k,z)$ is computed at $z=x$, where $x$ is any number such that $|x|>1$. Then continuing analytical function $J^*(k,x)$ along the intervals 
$(-\infty, -1)$, $(1,\infty)$ on the plane of complex variable $z$ with cut along the interval $[-1,1]$, we obtain
\begin{eqnarray}
J^*(k,z)=\frac{1}{\sqrt{z^2-1}}\left. arcsin\frac{zt-1}{z-t}\right|_{t=t_k}^{t=t_{k+1}}
=\frac{\pi h_k}{\sqrt{z^2-1}}F_k(z). \label{eq15}
\end{eqnarray}
For $J_1^*$
\begin{eqnarray}
J_1^*(k,z)=zJ^*(k,z)+g_k\pi h_k. \label{eq16}
\end{eqnarray}
Replacing $x$ into $z$ in \re{eq12} and using \re{eq14}-\re{eq16}, we obtain the coefficients $B_k$ of the QFs \re{eq4}

\section{Estimation of errors}

Let us introduce the following classes of functions:

\begin{enumerate}

\item $H^\a([-1,1],K)$ is a class function satisfying Holder condition on the interval
with the index $\a$ and constant $K$.

\item $C^{m,\a}[-1,1]=\Bigl\{f(t): \ \ f^{(m)} \in H^\a([-1,1],K_m)\Bigr\}$

\item $CC_{\triangle}^{m,\a}[-1,1]=\Bigl\{f(t): \ \ f(t) \in C[-1,1] \ \ \mbox{and} \ \
f(t) \in C^{m,\a}[t_k, t_{k+1}]  \Bigr\}$.

\item $W^r[-1,1]=\left\{f(t): f^{(r-1)}(t) \ \ \mbox{is absolutely continuous and } \ \
ess \sup\limits_{|t|\leq 1} |f^{(r)}|=M_r \right\} $.

\item $CW_{\triangle}^{r}[-1,1]=\Bigl\{f(t): \ \ f(t) \in C[-1,1] \ \ \mbox{and}  \ \
f(t) \in W^r[t_k,t_{k+1}] \Bigr\}$.

\end{enumerate}

\noindent Everywhere we use the notation $||f||_C=||f(t)||_{C[-1,1]}$ as a norm of the function. Note that $M_r=||f^{(r)}||_C$ for any $f(t) \in C^r[-1,1]$.

\noindent Now we prove the following theorems with respect to QFs \re{eq3} and \re{eq4}.
\begin{thm}
Let $f(t)$ be a function belonging to one of the classes of functions $W^1[-1,1]$, $CC_{\triangle}^{1,\a}$ or $CW_{\triangle}^{2}[-1,1]$.
Then for the errors of QFs \re{eq3} the estimations
$$
||R_N(f,x)||_C \leq L\frac{LnN}{N^\b},
$$
are true for all $x \in (-1,1)$, where $L$ and $\b$ are given in the Table 1.
\begin{table}[!ht]
\begin{center}
\caption{\small For QFs \re{eq3}} \label{t1}
\begin{tabular}{|c|c|c|}
  \hline
  Classes of functions & $\b$ & L \\
  \hline
  $W^1[-1,1]$ & 1 & $ \ds \frac{4\ga M_1}{\pi}\left(1+\frac{\pi\sqrt{2}}{2\ga lnN}\right)$ \\[3mm]

  $CC_{\triangle}^{1,\a}[-1,1]$ & $1+\a$ & $ \ds \frac{2\ga^{1+\a} K_1}{\pi}\left(1+\frac{\pi\sqrt{2}}{2\ga lnN}\right)$ \\[5mm]

  $CW_{\triangle}^{2}[-1,1]$ & 2 & $ \ds \frac{\ga^2 M_2}{\pi}\left(1+\frac{\pi\sqrt{2}}{\ga lnN}\right)$  \\[5mm]
  \hline
\end{tabular}
\end{center}
\end{table}
\end{thm}
{\sf Remark 1}: In the case of uniform grids, $\ga=2$.
\begin{thm}
Let $f(t)$ satisfy the conditions of Theorem 1. Then the errors of QFs \re{eq4} are
$$
\max\limits_{z}|R_N^*(f,z)| \leq L^*\frac{LnN}{N^\b}, \ \ L^*=\sqrt{L^2+\left(\frac{L_1}{lnN}\right)^2}
$$
where $L$, $L_1$ and $\b$ are given in the Table 2.
\begin{table}[!ht]
\begin{center}
\caption{\small For QFs \re{eq4}} \label{t2}
\begin{tabular}{|c|c|c|c|}
  \hline
  Classes of functions & $\b$ & L & $L_1$ \\
  \hline
  $W^1[-1,1]$ & 1 & $ \ds \frac{4\ga M_1}{\pi}\left(1+\frac{\pi\sqrt{2}}{2\ga lnN}\right)$ & $\ds \frac{M_1\ga}{2\pi}$ \\[5mm]

  $CC_{\triangle}^{1,\a}[-1,1]$ & $1+\a$ & $ \ds \frac{2\ga^{1+\a} K_1}{\pi}\left(1+\frac{\pi\sqrt{2}}{2\ga lnN}\right)$
  & $\ds \frac{K_1\ga^{1+\a}}{4\pi}$ \\[5mm]

  $CW_{\triangle}^{2}[-1,1]$ & 2 & $ \ds \frac{\ga^2 M_2}{\pi}\left(1+\frac{\pi\sqrt{2}}{\ga lnN}\right)$ & $\ds \frac{M_2\ga^2}{8\pi}$ \\[5mm]
  \hline
\end{tabular}
\end{center}
\end{table}

\end{thm}
Next theorem is again related to QFs \re{eq3} but in different classes of functions:
\begin{thm}
Let $f(t)$ be a function belonging to one of the classes of functions
$$
H^\a([-1,1],K), \ W^1[-1,1], \ CC_{\triangle}^{1,\a} \text{ or } CW_{\triangle}^{2}[-1,1].
$$
Then the error terms of QFs \re{eq3} satisfy the following estimations
$$
||R_N(f,x)||_C \leq L_2\frac{\ln N}{N^\b},
$$
for all $x \in (-1,1)$, where $L_2$ and $\b$ are given in the Table 3.
\begin{table}[!ht]
\begin{center}
\caption{\small For QFs \re{eq3}} \label{t3}
\begin{tabular}{|c|c|c|}
  \hline
  Classes of functions & $\b$ & $L_2$ \\
  \hline
  $H^\a([-1,1],K)$ & $\a$ & $ \ds \frac{2^{2-\a}\ga^\a K}{\pi}
  \left(1+\left(2+\frac{1}{\a}\right)\frac{2^{2-2\a}}{\ga^\a lnN}\right)$ \\[5mm]

  $W^1[-1,1]$ & 1 & $ \ds \frac{2\ga M_1}{\pi}\left(1+\frac{12\pi}{\ga lnN}\right)$ \\[5mm]

  $CC_{\triangle}^{1,\a}[-1,1]$ & $1+\a$ & $ \ds \frac{\ga^{1+\a} K_1}{\pi}\left(1+\frac{12\pi}{\ga lnN}\right)$ \\[5mm]

  $CW_{\triangle}^{2}[-1,1]$ & 2 & $ \ds \frac{\ga^2 M_2}{2\pi}\left(1+\frac{\pi\sqrt{24\pi}}{\ga lnN}\right)$  \\[5mm]
  \hline
\end{tabular}
\end{center}
\end{table}
\end{thm}
{\sf Remark 2}: Note that the main terms of $L_2$ in the Theorem 3 for the last three classes of functions
is twice less than the main terms of $L$ in the Theorem 1.

In estimation of the error of QFs \re{eq3} we use the idea of \cite{MS} (see also \cite{Ga2}) and the following Lemmas.

\setcounter{thm}{0}
\begin{lem}
Let $S_N(t)$ be linear spline \re{eq8} interpolating $f(t)$ on the grid $\triangle$, and let $t\in [t_k, t_{k+1}]$. Then for the estimate of error $r_N(f,t)=S_N(t)-f(t)$, we obtain
$$
||r_N(f,t)||_C=r_N^*(h_k),
$$
where $r_N^*(h_k)$ is defined by the Table \re{t4}
\begin{table}[!ht]
\begin{center}
\caption{\small Error of the linear spline \re{eq8} } \label{t4}
\begin{tabular}{|c|c|}
  \hline
  Classes of functions & $r_N^*(h_k)$ \\
  \hline

  $H^\a([-1,1],K)$ & $ \ds \frac{1}{2^{\a}}Kh_k^\a $ \\[5mm]

  $W^1[-1,1]$ & $ \ds \frac{1}{2}h_k$ \\[5mm]

  $CC_{\triangle}^{1,\a}[-1,1]$ & $ \ds \frac{1}{4}K_1h_k^{1+\a}$ \\[5mm]

  $CW_{\triangle}^{2}[-1,1]$ & $ \ds \frac{1}{8}h_k^2$  \\[5mm]
  \hline
\end{tabular}
\end{center}
\end{table}
\end{lem}
Lemma 1 is proved as Theorem 2.1 which is shown in \cite{MS}.
\begin{lem}
Let $S_N(t)$ be linear spline  defined by \re{eq8}. Then
$$
r_N(f,t) \in H^1([-1,1],\widetilde{K}), 
$$
where $\widetilde{K}$ is given in Table \re{t5}.
\begin{table}[!ht]
\begin{center}
\caption{\small Error of the linear spline \re{eq8} } \label{t5}
\begin{tabular}{|c|c|}
  \hline
  Classes of functions & $\widetilde{K}$ \\
  \hline

  $W^1[-1,1]$ & $ \ds 2M_1$ \\[5mm]

  $CC_{\triangle}^{1,\a}[-1,1]$ & $ \ds K_1h^\a_k$ \\[5mm]

  $CW_{\triangle}^{2}[-1,1]$ & $ \ds M_2h_k$  \\[5mm]
  \hline
\end{tabular}
\end{center}
\end{table}
\end{lem}
\noindent {\bf Proof of the Lemma 2}. Consider three cases:
$$
(a) t,t' \in [t_k, t_{k+1}], \ \ (b) |t-t'|\geq h_k
\ \ (c) \tau \in [t_{k-1}, t_k], \ \  t' \in [t_{k}, t_{k+1}], \ \ |t-t'| \leq h_k
$$
\textbf{I}. Let $f(t) \in W^1[-1,1]$. Then in the case $(a)$, using representation \re{eq8} we have
\begin{eqnarray}
&& \ds \Bigl|r_N(f;t)-r_N(f;t')\Bigr|=\frac{1}{h_k}\Bigl|(t-t')[f(t_{k+1})-f(t_k)]-
h_k[f(t)-f(t')]\Bigr| \nn \\
&& \qquad \qquad  \ds =\frac{1}{h_{k}}\left|(t-t')\int\limits_{t_k}^{t_{k+1}}f'(s)ds -
h_k \int\limits_{t'}^{t}f'(s)ds\right| \leq 2M_1|t-t'|.  \nn
\end{eqnarray}
In the case $(b)$, in accordance with Lemma 1, we get
\begin{eqnarray}
&& \ds \Bigl|r_N(f;t)-r_N(f;t')\Bigr| \leq \Bigl|r_N(f;t)\Bigr| + \Bigl|r_N(f;t')\Bigr|
\leq h_kM_1 \leq M_1 |t-t'|.  \nn
\end{eqnarray}
Using the case $(a)$, in the case $(c)$ we obtain
\begin{eqnarray}
&& \ds \Bigl|r_N(f;t)-r_N(f;t')\Bigr| \leq \Bigl|r_N(f;t)-r_N(f;t_k)\Bigr| + \Bigl|r_N(f;t_k)-r_N(f;t')\Bigr| \nn \\
&& \qquad \qquad \qquad \qquad \quad \ds \leq 2M_1|t_k-t|+2M_1|t-t'|= 2M_1 |t-t'|.  \nn
\end{eqnarray}
\textbf{II}. Now let  $f(t) \in CC_{\triangle}^{1,\a}[-1,1]$. In the case $(a)$
\begin{eqnarray}
&& \ds \Bigl|r_N(f;t)-r_N(f;t')\Bigr| \leq \frac{1}{h_k}
\Bigl|(t-t')[f(t_{k+1})-f(t_k))-h_k[f(t)-f(t')]\Bigr|  \nn \\
&& \qquad \qquad \qquad \qquad \quad \ds = \frac{1}{h_k}
\Bigl|(t-t')(t_{k+1}-t_k)f'(\theta_1)-h_k(t-t')f'(\theta_2)\Bigr|  \nn \\
&& \qquad \qquad \qquad \qquad \quad \ds =
\Bigl|(t-t')\Bigr|\Bigl|f'(\theta_1)-f'(\theta_2)\Bigr|
\leq K_1 |t-t'||\theta_1-\theta_2|^\a \nn \\
&& \qquad \qquad \qquad \qquad \quad \ds \leq K_1h_k^\a|t-t'|. \nn
\end{eqnarray}
In the case $(b)$, due to Lemma 1, we have
\begin{eqnarray}
&& \ds \Bigl|r_N(f;t)-r_N(f;t')\Bigr| \leq K_1 h_k^\a
|t-t'|.  \nn
\end{eqnarray}
It is obvious in the case $(c)$ that
\begin{eqnarray}
&& \ds \Bigl|r_N(f;t)-r_N(f;t')\Bigr| \leq K_1 h_k^\a
|t-t'|.  \nn
\end{eqnarray}
\textbf{III}. Let $f(t) \in CW_{\triangle}^{2}[-1,1]$. In case $(a)$, we have 
\begin{eqnarray}
&& \ds \Bigl|r_N(f;t)-r_N(f;t')\Bigr| = \Bigl|t-t'\Bigr|
\Bigl|f'(\theta_1)-f(\theta_2)\Bigr|  \nn \\
&& \qquad \qquad \qquad \qquad \quad \ds
= \Bigl|t-t'\Bigr|\left|\int\limits_{\theta_1}^{\theta_2}f''(s)ds\right|
\leq M_2h_k|t-t'|. \nn
\end{eqnarray}
The cases $(b)$ and $(c)$ are proved in a similar way as the case $(a)$. So that the proof of the Lemma 2 follows
from the above obtained errors.

\noindent Now it is easy to prove the following lemma.
\begin{lem}
Let $S_N(t)$ be linear spline defined by \re{eq8} and $f(t) \in H^\a([-1,1],K)$. Then
$$
r_N(f,t) \in H^\a([-1,1],2^{2-\a}K).
$$
\end{lem}

\noindent {\bf Prove of the Theorem 1}. Since
$$
\int\limits_{-1}^{1}\frac{dt}{\sqrt{1-x^2}(t-x)}=0,
$$
the reminder term of QFs \re{eq3} can be represented as
\begin{eqnarray} \label{eq17}
&& \ds R_N(f,x) = \frac{\sqrt{1-x^2}}{\pi}\int\limits_{-1}^{1}\frac{r_N(f,t)-r_N(f,x)}{\sqrt{1-x^2}(t-x)}dt.
\end{eqnarray}
For definiteness, let us prove the Theorem 1 in case $0\leq x \leq 1$ (the case $-1 \leq x \leq 0$ is considered analogically).
Fixing the number $0 < \de_N < \frac{1}{2}$ and dividing the integral in \re{eq17} into three parts to yield
\begin{eqnarray} \label{eq18}
&& \ds R_N(f,x) = \frac{\sqrt{1-x^2}}{\pi}\left(\int\limits_{-1}^{x-\de_N} + \int\limits_{x-\de_N}^{x+\de_N} +
\int\limits_{x+\de_N}^{1} \right) \frac{r_N(f,t)-r_N(f,x)}{\sqrt{1-x^2}(t-x)}dt \nn \\
&& \ds \qquad \qquad = \frac{\sqrt{1-x^2}}{\pi}\left( J_1 + J_2 + J_3 \right).
\end{eqnarray}
First assume that $\de_N < 1-x$. Then due to \re{eq10}, for $J_1$ we have
\begin{eqnarray}
&& \ds \Bigl| J_1 \Bigr| = 2 ||r_N(f,x)||_C\left|\int\limits_{-1}^{x-\de_N}\frac{dt}{\sqrt{1-x^2}(t-x)}\right| \nn \\
&& \ds \qquad = \frac{2}{\sqrt{1-x^2}} ||r_N(f,x)||_C ln\left|\frac{t\sqrt{1-x^2} - x\sqrt{1-t^2}}{\sqrt{1-x^2} +\sqrt{1-x^t}}\right|_{t=x-\de_N}. \nn
\end{eqnarray}
It is not hard to show that
\begin{eqnarray}
&& \vp_1(x,\de_N) = \left. \frac{t\sqrt{1-x^2}-x\sqrt{1-t^2}}{\sqrt{1-x^2} +\sqrt{1-x^t}}\right|_{t=x-\de_N} \nn \\
&& \ds \qquad \qquad = x + \frac{1}{\de_N} \left( 1-x^2 -\sqrt{1-x^2}\sqrt{1-(x-\de_N)^2}\right). \nn
\end{eqnarray}
This is a function of $x$ which strictly decreases on $\left[0, \frac{\de_N}{2}\right]$ and strictly increases on $\left[\frac{\de_N}{2}, 1\right]$, and
$$
\vp_1(0,\de_N) = \frac{\de_N}{1+\sqrt{1-\de_N^2}}, \ \ \vp_1(\de_N/2,\de_N) = \frac{\de_N}{2}, \ \ \vp_1(1,\de_N) = 1.
$$
Hence,
\begin{eqnarray} \label{eq19}
|J_1| \leq \frac{2}{\sqrt{1-x^2}}||r_N(f,x)||_C ln\frac{2}{\de_N}.
\end{eqnarray}
For $J_2$, we use Lemma 2
\begin{eqnarray}
&& |J_2| \leq \tilde{K}\int\limits_{x-\de_N}^{x+\de_N} \frac{dt}{\sqrt{1-t^2}} = \tilde{K}
[\arcsin(x+\de_N) - \arcsin(x-\de_N)]. \nn
\end{eqnarray}
Let
$$
\vp_2(x,\de_N) =\arcsin(x+\de_N) - arcsin(x-\de_N).
$$
Since $0 \le x \leq 1-\de_N $ and by assumption $\de_N \leq 1-x$, derivative of $\vp_2(x,\de_N)$ is positive and
$$
\vp_2(x,\de_N) \leq \vp_2(1-\de_N, \de_N) = \arcsin 2\sqrt{\de_N(1-\de_N)}.
$$
From this and the known inequality $\arcsin \a \leq \frac{\pi}{2}\a, (0\leq \a \leq \frac{\pi}{2})$ it follows that
$$
\vp_2(x,\de_N) \leq \pi \sqrt{\de_N}.
$$
Hence
\begin{eqnarray} \label{eq20}
&& |J_2| \leq \tilde{K}\pi \sqrt{\de_N}
\end{eqnarray}
For $J_3$, we have
\begin{eqnarray}
&& |J_3| = 2||r_N(f,x)||_C \left|\int\limits_{x+\de_N}^{1}\frac{dt}{\sqrt{1-x^2}(t-x)} \right| \nn \\
&& \ds \qquad \qquad = \frac{2}{\sqrt{1-x^2}} ||r_N(f,x)||_C (-1) ln\left|\frac{t\sqrt{1-x^2} - x\sqrt{1-t^2}}
{\sqrt{1-x^2} + \sqrt{1-t^2}}\right|_{t=x+\de_N}. \nn
\end{eqnarray}
We may show that the function
\begin{eqnarray}
\vp_3(x,\de_N)= \left. \frac{t\sqrt{1-x^2} - x\sqrt{1-t^2}}{\sqrt{1-x^2} + \sqrt{1-t^2}} \right|_{t=x+\de_N}
= -\vp_1(x, -\de_N), \nn
\end{eqnarray}
strictly  increases on $[0, 1-\de_N]$ and strictly decreases from $\frac{\de_N}{1+\sqrt{1-\de_N^2}}$ to 1. So that
\begin{eqnarray} \label{eq21}
&& |J_3| \leq \frac{2}{\sqrt{1-x^2}}||r_N(f,x)||_C \ln\frac{2}{\de_N}.
\end{eqnarray}
It follows from the errors of \re{eq19}- \re{eq21} and \re{eq18} that
\begin{eqnarray} \label{eq22}
&& ||R_N(f,x)||_C \leq \frac{4}{\pi}||r_N(f,x)||_C \ln\frac{2}{\de_N} +\tilde{K} \sqrt{\de_N}.
\end{eqnarray}
Now consider the case $\de_N > 1-x$. Write
\begin{eqnarray} \label{eq23}
&& \ds R_N(f,x) = \frac{\sqrt{1-x^2}}{\pi}\left(\int\limits_{-1}^{x-\de_N} + \int\limits_{x-\de_N}^{1} \right)
\frac{r_N(f,t)-r_N(f,x)}{\sqrt{1-x^2}(t-x)}dt \nn \\
&& \ds \qquad \qquad = \frac{\sqrt{1-x^2}}{\pi}\left( J_1^* + J_2^* \right)
\end{eqnarray}
Integral $J_1^*$ is estimated as $J_1$. Due to Lemma 2
\begin{eqnarray}
&& |J_2^*| \leq \tilde{K}\int\limits_{x-\de_N}^{1}\frac{dt}{\sqrt{1-t^2}} = \tilde{K} \arcsin\sqrt{1-(x-\de_N)^2}. \nn
\end{eqnarray}
Since $ 0< x \leq 1, \ \ \de_N > 1 -x $ and due to the inequality
$$
1-(x-\de_N)^2 = (1-x +\de_N)(1+x-\de_N) < 4\de_N.
$$
we obtain
$$
|J_2^*| \leq \tilde{K} \arcsin2\sqrt{\de_N} \leq \tilde{K} \pi \sqrt{\de_N}.
$$
Substituting the errors of $J_1^*$ and $J_2^*$ into \re{eq23}, we arrive at estimation \re{eq22}.

\noindent In order to determine the errors of estimation for every classes of functions in Theorem 1, we use the results of Lemma 1 and 2
and set $\de_N=\frac{2}{N^2}$. Viz: \\

\noindent \textbf{I}. Let $f(t) \in W^1[-1,1]$. Then
\begin{eqnarray}
&& \ds ||R_N(f,x)||_C \leq \frac{2}{\pi}M_1 h \ln\frac{2}{\de_N} + 2M_1\sqrt{\de_N} \nn \\
&& \ds \qquad \qquad \qquad = \frac{4M_1\ga}{\pi}\left(1+\frac{\pi \sqrt{2}}{2\ga \ln N} \right)\frac{\ln N}{N}. \nn
\end{eqnarray}
\noindent {\bf II}. If $f(t) \in CC_\triangle^{1,\a}[-1,1]$, then
\begin{eqnarray}
&& \ds ||R_N(f,x)||_C \leq \frac{K_1}{\pi}h^{1+\a}\ln\frac{2}{\de_N} + 2K_1h^\a\sqrt{\de_N} \nn \\
&& \ds \qquad \qquad \qquad = \frac{2K_1\ga^{1+\a}}{\pi}\left( 1 + \frac{\pi \sqrt{2}}{2\ga \ln N} \right)\frac{\ln N}{N^{1+\a}}. \nn
\end{eqnarray}
\noindent {\bf III}. If $f(t) \in CW_\triangle^{2}[-1,1]$, then
\begin{eqnarray}
&& \ds ||R_N(f,x)||_C \leq \frac{M_2}{2\pi}\frac{\ga^2}{N^2}\ln\frac{2}{\de_N} + \frac{2M_2 \ga}{N}\sqrt{\de_N} \nn \\
&& \ds \qquad \qquad \qquad = \frac{M_2\ga^{2}}{\pi}\left( 1 + \frac{\pi \sqrt{2}}{\ga \ln N} \right)\frac{\ln N}{N^{2}}. \nn
\end{eqnarray}
Theorem 1 is proved.

\textbf{Proof of the Theorem 2}  is carried out by the famous scheme of the formula Sokhotskii-Plemergh (see \cite{Mu}), principle maximum module for analytical function and results of Theorem 1 and Lemma 1.

\textbf{Proof of the Theorem 3}. Let the remainder term of QFs \re{eq3} be divided into three  parts
\begin{eqnarray} \label{eq18}
&& \ds R_N(f,x) = \frac{\sqrt{1-x^2}}{\pi}\left(\int\limits_{-1}^{x-\de_N} + \int\limits_{x-\de_N}^{x+\de_N} +
\int\limits_{x+\de_N}^{1} \right) \frac{r_N(f,t)-r_N(f,x)}{\sqrt{1-x^2}(t-x)}dt \nn \\
&& \ds \qquad \qquad = \frac{\sqrt{1-x^2}}{\pi}\left( \widetilde{J}_1 + \widetilde{J}_2 + \widetilde{J}_3 \right).
\end{eqnarray}
In the proof of Theorem 1, we have already seen that the case $\de_N < 1-x$, is adequate for the estimations of $\widetilde{J}_1$
and $\widetilde{J}_3$ i.e.
\begin{eqnarray} \label{eq25}
&& |\widetilde{J}_1| + |\widetilde{J}_3| \leq \frac{4}{\pi}\sqrt{1-x^2}||r_N(f,x)||_C \ln \frac{2}{\de_N}.
\end{eqnarray}
For the estimation of $\widetilde{J}_2$, we consider the function
\begin{eqnarray} \label{eq26}
&& T(x,\ve, \s) = \sqrt{1-x^2}\int\limits_{x-\ve}^{x+ve}\frac{|t-x|^{\a-1}}{\sqrt{1-t^2}}dt,
\end{eqnarray}
where $0 < \s \leq 1, \ \  1 - x \geq \ve, \ \ 0<\ve <\frac{1}{2}$.

\noindent It is obvious that
\begin{eqnarray} \label{eq27}
&& \ds T(x,\ve,\s) \leq 2\ve^\s \int\limits_0^1\frac{\sqrt{1-x^2}y^{\s-1}}{\sqrt{1-(x+\ve y)^2}}dy \nn \\
&& \ds \qquad \qquad = 2\ve^\s \left[\int\limits_0^{1/2}\frac{\sqrt{1-x^2}y^{\s-1}}{\sqrt{1-(x+\ve y)^2}}dy +
\int\limits_{1/2}^{1}\frac{\sqrt{1-x^2}y^{\s-1}}{\sqrt{1-(x+\ve y)^2}}dy\right] \nn \\
&& \ds \qquad \qquad = 2\ve^\s[T_1(x,\ve,\s) + T_2(x,\ve,\s)].
\end{eqnarray}
Let $k_0 = [1/\ve]$, since $1-x \geq \ve $, then for some $k, 1\leq k \leq k_0$ the inequality 
$$ k\ve \leq 1-x < (k+1)\ve, $$
takes place. From this and $\ds 0 < y \leq \frac{1}{2} $ it follows that
$$
\frac{1-x^2}{1-(x+\ve y)^2} = \frac{(1-x)(1+x)}{(1-x-\ve y)(1+ x + \ve y)} \leq \frac{1-x}{1-x-\ve/2} \leq \frac{k+1}{k-1/2} \leq 4,
$$
for all $k\geq 1$. Hence
\begin{eqnarray} \label{eq28}
T_1(x,\ve,\s) = 2\int\limits_0^{1/2}y^{\a-1}dy = \frac{2^{1-\a}}{\s}.
\end{eqnarray}
Furthermore
\begin{eqnarray}
T_2(x,\ve,\s) = \frac{1}{2^{\s-1}}\int\limits_{1/2}^1 \frac{\sqrt{1-x^2}}{\sqrt{1-(x+\ve y)}}dy. \nn
\end{eqnarray}
Let  $\ve \leq 1-x < 2\ve$. Then
$$
\frac{1-x^2}{1-(x+\ve y)^2} \leq \frac{1-x}{1-x-\frac{\ve}{2}} \leq \frac{2}{1-y},
$$
and therefore
\begin{eqnarray} \label{eq29}
T_2(x,\ve,\s) = 2^{1-\s}\sqrt{2}\int\limits_{1/2}^1 (1-y)^{1/2}dy = 2^{2-\s}.
\end{eqnarray}
If for some $k \leq 2$ the inequality
$$
k\ve \leq 1-x \leq (k+1)\ve
$$
takes place, then
$$
\frac{1-x^2}{1-(x+\ve y)^2} \leq \frac{1-x}{1-x- \ve y} \leq \frac{k+1}{k-1} \leq 3,
$$
therefore
\begin{eqnarray} \label{eq30}
T_2(x,\ve,\s) \leq 2^{1-\s}\sqrt{3} \cdot \frac{1}{2} < 2^{2-\s}.
\end{eqnarray}
From \re{eq27}-\re{eq30} it follows that
\begin{eqnarray} \label{eq31}
T(x,\ve,\s) \leq 2^{1-\s}\left(2+\frac{1}{\s}\right) \ve^\s.
\end{eqnarray}
Now for $\widetilde{J}_2$ we have
\begin{eqnarray} \label{eq32}
\widetilde{J}_2 \leq \frac{K}{\pi} \sqrt{1-x^2}\int\limits_{x-\de_N}^{x+\de_N}\frac{|t-x|^{\s-1}}{\sqrt{1-t^2}} = \frac{K}{\pi} T(x,\ve,\s),
\end{eqnarray}
where $\s=\a$ for the class $H^\a([-1,1],K)$ and $\s = 1$ for rest classes of functions. Assuming $\s = \frac{2}{N}$ and from the Lemmas 2, 3 and inequalities \re{eq25} and \re{eq31}-\re{eq32} we get the assertion of the Theorem 3.


\begin{thebibliography}{EIBLM}

\bibitem{BL}
S. Belotserkovskii, I. Lifanov (1985) {\em Numerical methods in
singular integral equations. Moscow: Nauka, 254 p. (in Russian).}

\bibitem{Ga1} Gabdulkhaev B.G. Finite-dimentional approximation of singular integrals and direct methods
for special integrals and integra-differential equations. Mathematical analysis. Resume science and technology.
VINITI AN SSSR, V.18, pp. 251-307.

\bibitem{Ga2} Gabdulkhaev B.G. Optimal approximation for the solution of linear problems.
Kazan: Kazan University press, 1980. 231 p.

\bibitem{MS}
Makavoz Yu. I., Sheshko M.A. On a error of estimation of quadrature formulas for singular integrals.
Isv. AN SSSR. Ser. Phiz-Math. Nauk. 1977. V.6, pp. 36-41.

\bibitem{Mu}
Muskhelishvili N.I. (1953) {\em Singular Integral equations.
Gostekhizda (1946). M: Fizmatgiz, 1962. 512p. }

\bibitem{Pi}
Pikhtiev G.H. Accurate methods for evaluation of Cauchy type singular integrals. Novosibirsk: Science, 1980.

\bibitem{ZKM} U.S.Zavyalov, B.I.Kvasov, B.I.Miroshnichenko. Methods of spline functions,
Nauka,Moskov,1980

\end{thebibliography}
\end{document}